\documentclass[12pt]{article}
\usepackage[mathscr]{eucal}
\usepackage{amsbsy}
\usepackage{amssymb}
\usepackage{amsmath}
\usepackage{amsthm}
\usepackage{graphics}
\usepackage{epsfig}

\usepackage{latexsym}

\usepackage[left=2cm, top=2cm, right=1.5cm, bottom=2cm] {geometry}

\newcommand{\PP}{\mathbf{P}}

\def\phi{\varphi}
\def\bbn{{\mathbb N}}
\def\bbr{{\mathbb R}}

\def\bz{{\bf Z}}

\def\bz{{\mathbb Z}} %Keistas Z apribrezimas

\def\Bl1{{\bf 1}}
\def\B2{{\bf 2}}
\def\B0{{\bf 0}}

\def\a{\alpha}
\def\b{\beta}
\def\d{\delta}

\def\e{\varepsilon}
\def\g{\gamma}

\def\l{\lambda}

\def\=A8{\"o}

\def\pc{\stackrel{P}{\longrightarrow}}

\def\fdd{\stackrel{f.d.d.}{\longrightarrow}}

\newcommand{\beq}{\begin{equation}}
\newcommand{\eeq}{\end{equation}}
\newcommand\beqn{\begin{displaymath}}  % no number
\newcommand\eeqn{\end{displaymath}}

\newcommand{\dx}{{\rm d}x }

\newcommand{\ind}[1]{\mathbbm{1}_{#1}}

\newcommand{\halmos}{\vspace{3mm} \hfill \mbox{$\Box$}\\[2mm]}

\theoremstyle{plain}
\newtheorem{teo}{Theorem}

\newtheorem{prop}[teo]{Proposition}
\theoremstyle{definition}

\newtheorem{remark}[teo]{Remark}

\usepackage{bbm}
\begin{document}

\title{A note on linear processes with tapered innovations
\footnotemark[0]\footnotetext[0]{ \textit{Short title:}
linear processes with tapered innovations }
\footnotemark[0]\footnotetext[0]{%
\textit{MSC 2010 subject classifications}. Primary 60G99, secondary
60G22, 60F17 .} \footnotemark[0]
\footnotemark[0]\footnotetext[0]{ \textit{Corresponding author:}
Vygantas Paulauskas, Department of Mathematics and
 Informatics, Vilnius university, Naugarduko 24, Vilnius 03225, Lithuania,
 e-mail:vygantas.paulauskas@mif.vu.lt}}

\author{ Vygantas Paulauskas$^{\text{\small 1}}$ \\
{\small $^{\text{1}}$ Vilnius University, Department of Mathematics
and
 Informatics,}\\}

\date{April 20, 2019}

%\begin{document}

\maketitle

\begin{abstract}

In the paper we consider the partial sum process $\sum_{k=1}^{[nt]}X_k^{(n)}$, where $\{X_k^{(n)}, \ k\in \bz\},\ n\ge 1,$ is a series of linear processes with innovations having heavy-tailed tapered distributions with tapering parameter $b_n$ depending on $n$. It is shown that, depending on the properties of a filter of a linear process under consideration and  on the parameter $b_n$ defining if the tapering is hard or soft, the limit process for such partial sum process can be fractional Brownian motion or linear fractional stable motion.

\end{abstract}
%\vfill
%\eject
\section{Introduction and formulation of the results}

Recently the so-called trawl processes were introduced, see \cite{Barndorff1}, \cite{Barndorff2}, \cite{Doukhan}.
In the first two cited papers more general processes of continuous time, giving trawl processes as special case, were considered, in the third one the trawl processes of discrete time
were introduced. These stationary random processes are defined as follows. Suppose we have  a  random process $\g=\{\g(u), u\in \bbr\}$, tending to zero in probability as $u\to 0$, and deterministic sequence $\{a_j\in \bbr, j\in \bbn\}$,  $ \lim_{j\to \infty} a_j=0$. Random process $\g$ is called the seed process and $\{a_j, j\in \bbn\}$  is called the trawl. Then the trawl process, corresponding to the seed process $\g$ and the trawl $\{a_j, j\in \bbn\}$, is defined as
\begin{equation}\label{trawl}
X_k=\sum_{j=0}^\infty \g_{k-j}(a_j), \ k\in \bz,
\end{equation}
where $\g_k=\{\g_k(u), u\in \bbr\}$ are independent and identically distributed (i.i.d) copies of $\g$. Clearly, some additional conditions on the seed process and the trawl must be required in order to ensure the a.s. convergence of the series (\ref{trawl}).  Taking different seed processes (Brownian motion, Poisson or Bernoulli  processes) we get stationary sequences with different properties of dependence and different limit behaviour for partial sum processes $\sum_{j=1}^{[nt]}X_j$. The usual linear processes
\begin{equation}\label{linpr1}
X_k=\sum_{j=0}^\infty a_j\xi_{k-j}, \ k\in \bz,
\end{equation}
are obtained taking the most simple seed process  $\g=\{\xi u, u\in \bbr\}$ with some random variable $\xi$ (in \cite{Doukhan} it is called the random line seed process), therefore, $\xi_i, i\in \bz,$ are i.i.d. copies of $\xi$.
For linear processes we shall use the terminology which is traditional for this class of random processes: a sequence $\{a_j, j\in \bbn\}$ and a sequence of i.i.d. random variables $\xi_i, i\in \bz,$ will be called the filter and innovations, respectively,   of a linear process $\{X_k, \ k\in \bz\}$. Depending on the moments of innovations we can model stationary sequences $\{X_k, \ k\in \bz\}$ with finite or infinite variance, while properties of a filter are responsible for the dependence structure and memory properties (with respect to the summation operation) of this sequence, thus we can have short-range or long-range dependencies (the case of sequences with finite variance) and four types of memory -positive, zero, negative, and strongly negative, see \cite{Paul20}, where this terminology is suggested, instead of long, short and negative memories. Limit theorems for partial sum processes $\sum_{j=1}^{[nt]}X_j$, formed by linear processes, are well-studied object, documented in many papers and several monographs. Instead of simple random line seed process (which gives a linear process) taking more complicated seed processes we obtain more complicated stationary sequences. Limit theorems
for partial sum processes  formed by such stationary sequences are studied in \cite{Doukhan}. It is shown  (see Theorems 1 and 2 therein) that the limit processes can be the fractional Brownian motion (fBm) or an $\a$-stable L{\' e}vy process.

In this note we want to demonstrate that similar behavior of partial sum processes, formed by trawl processes, can be observed in the context of linear processes (the most simple example of the trawl process), but using the so-called tapered innovations or/and tapered filter.
Let us introduce the notion of the tapered distribution.  For more than a half of a century it was known that many natural hazards follow a power-law distribution. Among such hazards earthquakes are the most well-known hazard that can be described by power-law, but there are many papers which demonstrate that other natural hazards, such as rock falls, landslides, riverine floods, tsunami, wildfire exhibit power-law behavior, see, for example, \cite{Geist}. The power-law distribution, which we are speaking about, nowadays is called Pareto distribution, but for many years in seismology it was called G-R law (from the names of the authors of \cite{Gutenberg}). The Pareto probability density and distribution functions are given by the following expressions
\begin{equation}\label{paretodf}
f(x)= \a x_0^\a x^{-\a-1}, \quad F(x)=1-\left (\frac{x_0}{x}\right )^\a, \ x\ge x_0,
\end{equation}
where $\a>0, \ x_0>0$ are two parameters of Pareto distribution. For example, in seismology the scalar value of seismic moment of an earthquake follows Pareto distribution with some $x_0$ (minimal seismic moment, above which it is considered that seismic moment follows power-law) and exponent $\a$. Theoretical reasoning, based on the theory of branching processes, predicts that this exponent is universal constant equal $1/2$, although estimations from real data usually give a little bit bigger value, see discussion in \cite{Kagan}. Another important fact was that empirical data of earthquakes showed that while seismic moment follows power-law in quite big range of values, the largest values of samples demonstrate much lighter tails, see Figure 1 in \cite{Kagan}. Therefore it was suggested to apply the exponential taper to Pareto distribution, obtaining the following density and distribution functions
\begin{equation}\label{taperpd}
f(x, b)= (\a +xb^{-1})x_0^\a x^{-\a-1}\exp\left \{\frac{x_0-x}{b} \right \}, \ x\ge x_0,
\end{equation}
\begin{equation}\label{taperdf}
  F(x, b)=\left \{\begin{array}{ll}
  0,               & x< x_0, \\
  1-\left (\frac{x_0}{x}\right )^\a\exp\left \{\frac{x_0-x}{b} \right \}, & \ x\ge x_0,
  \end{array} \right.
\end{equation}
where $b>x_0$ is some parameter (in seismology it is called the corner moment), demonstrating the point where power-law is changing to the exponential one. Tapered Pareto (TP) distribution has moments of all orders, while Pareto distribution, which is obtained from (\ref{taperdf}) by letting $b\to \infty$, has moments of the order only up to $\a$.

We shall consider linear processes with innovations having the TP distribution, only for simplicity of calculations we shall consider random variables having  probability density function
\begin{equation}\label{taperpd1}
f_b(x)= \left \{\begin{array}{ll}
                0,     & {\rm for} \ x<1,\\
               \a x^{-\a-1}, & {\rm for} \ 1\le x\le b, \\
               b^{-\a}\exp (-x+b), & {\rm for} \  x> b.
               \end{array} \right.
\end{equation}
Such choice is made for the reason that parameter $x_0$, present in (\ref{taperpd}) , also the continuity of the density in (\ref{taperpd}) (density (\ref{taperpd1}) is discontinuous) are unimportant in our context.

Let  $\zeta$ stand for a random variable with the density (\ref{taperpd1}) and   $\xi=\zeta-E\zeta.$
Taking i.i.d. copies $\xi_k, \ k\in \bz$, of $\xi$ as innovations in (\ref{linpr1})  and assuming  the condition
\begin{equation}\label{cond}
\sum_{j=0}^\infty |a_j|^2 <\infty
\end{equation}
we get a family of linear processes, depending on the filter and two parameters $\a$ and $b$:
%\begin{equation}\label{linpr2}
$X_k=X_k(\a, b):=\sum_{j=0}^\infty a_j\xi_{k-j}, \ k\in \bz.$
%\end{equation}
Such  processes we call linear processes with tapered innovations. If we form the partial sum process
\begin{equation}\label{linpr5}
S_n(t)=\sum_{j=1}^{[nt]}X_j, \ 0\le t\le 1,
\end{equation}
with a fixed $b$, then we are in the well investigated situation of linear processes with finite variance, and, having regular behavior of the filter $\{a_j, j\in \bbn\}$, we get  the   fBm as a limit for the appropriately normalized $S_n(t)$. But the situation becomes different if we take the parameter $b$ depending on $n$ and growing unboundedly as $n\to \infty$. Namely, now we consider
\begin{equation}\label{linpr3}
{\bar S}_n(t)=\sum_{j=1}^{[nt]}X_j^{(n)}, \ 0\le t\le 1,
\end{equation}
where $X_j^{(n)}=X_j(\a, b_n)$. The limit behavior of ${\bar S}_n(t)$ is more complicated and depends on  three factors: the decay of a filter $\{a_i\}$, the growth of the tapering parameter $b_n$, and $\a$.

Although in the paper we shall not consider another possibility, let us mention that
the tapering procedure can be applied to a filter of a linear process, too. During last decade there were several papers devoted to this approach, we refer the reader to \cite{Sabzikar} and references in this paper. One can take the filter $c_j(\l):=a_j \exp (-j/\l),$ depending of the parameter $\l>1$,
and can consider linear processes with tapered filter (here it is necessary to note, that in the above cited  paper the term "tempered"is used instead of "tapered")
\begin{equation}\label{linpr4}
W_k=W_k(\l):=\sum_{j=0}^\infty c_j(\l)\eta_{k-j}, \ k\in \bz,
\end{equation}
and with some sequence of innovations $\{\eta_i, \ i\in \bz\}.$
As in the case of tapered innovations, we take     $W_k^{(n)}=W_k(\l_n)$ with parameter $\l_n \to \infty$  and  form the partial sum process:
$$
U_n(t)=\sum_{j=1}^{[nt]}W_j^{(n)}, \ 0\le t\le 1.
$$
 The limit behavior of $U_n$ will depend on the innovations $\{\eta_i\}$,  filter $\{a_i\}$, and the growth of the tapering parameter $\l_n$. Even we can consider partial sum processes formed by linear processes with tapered filter and innovations, if
we put in (\ref{linpr4}) $\xi_k$ with tapering parameter $b_n$ instead of $\eta_k$ , thus we get  partial sum process depending on two sequences of tapering parameters $b_n$ and $\l_n.$ In this paper we restrict ourselves only with linear processes with tapered innovations. We shall show that, depending on the growth of $b_n$ and normalizing sequence, the limit process for ${\bar S}_n(t)$ from (\ref{linpr3}) can be fBm or linear fractional stable motion.

 Let us note that during last decade there were several papers where truncated or tapered heavy-tailed distributions were considered. Truncation can be considered as a special case of tapering when we do not leave any mass over truncation level, while by tapering we only change (from the point $b$) the heavy tail of a distribution by the light tail of another distribution. First of all, it is necessary to mention the paper \cite{Chakrabarty}, where tapered (although in the paper the term "truncated" is used) heavy-tailed $\bbr^d$-valued random vectors were considered. In order not to introduce multivariate notation we take a particular case $d=1$, then it is possible to say that the following scheme of triangular array is considered in this paper:
 \begin{equation}\label{taper-trunk}
 \kappa_i^{(n)}=\theta_i\ind{[|\theta_i|<b_n]}+\frac{\theta_i}{|\theta_i|}(b_n+R_i)\ind{[|\theta_i|\ge b_n]},
 \end{equation}
 where $\{\theta_i \}$ are i.i.d. random variables with regularly varying tails and with the tail exponent $0<a<2$, and $\{R_i\}$ is a sequence of i.i.d. non-negative random variables , independent of the sequence $\{\theta_i \}$. Since in \cite{Chakrabarty} random variables $\kappa_i$ are called truncated, the
 following terminology is introduced:  a sequence of truncation levels $\{b_n\}$ is called soft truncation or hard truncation for a sequence $\{\theta_i, \ i\in N,\}$,  if
\begin{equation}\label{trunc}
\lim_{n\to \infty} n\PP\{|\theta_1|>b_n\}
\end{equation}
 is equal to $0$ or $\infty,$ respectively. In the paper \cite{Chakrabarty} asymptotic behavior of sums $\sum_{i=1}^n \kappa_i^{(n)} $ is investigated, and
 the main result  says that in the case of soft truncation sums of heavy-tailed truncated random variables behaves like sums of heavy-tailed ones, while in the case of hard truncation they behave like sums of light-tailed random variables.

  Since we think that it is more natural to call random variables $\{\kappa_i^{(n)} \} $ tapered, we shall use the terms soft and hard tapering. Truncation in (\ref{taper-trunk}) is obtained by setting $R_1$ equal to zero with probability $1$ and attributing all mass of a distribution of $\theta_1$ over the area $\{x\in \bbr: |x|>b_n\}$ to the points $\pm b_n$. Here it is worth to mention that in the literature one can find a little bit different truncation. For example, in \cite{Aban} the truncated Pareto distribution is defined as  a random variable with the density
  \begin{equation}\label{trunkPareto}
g_b(x)= \left \{\begin{array}{ll}
                \a x_0^\a x^{-\a-1}(1-(x_0/b))^{-1}, & {\rm for} \ 0<x_0\le x\le b<\infty, \\
                0,     & {\rm elsewhere}.
               \end{array} \right.
\end{equation}
In this truncation procedure the mass of Pareto distribution over point $b$ is distributed over interval $[x_0, b]$ by the same power-law. It would be possible to consider in our scheme such truncated Pareto distributions instead of TP distributions, defined by (\ref{taperpd1}), but the calculations would be more complicated, since the densities of Pareto and  truncated Pareto from (\ref{trunkPareto}) distributions differ on all interval $[x_0, \infty).$

In our case the initial sequence $\{\theta_i\}, \ i\in N,$ is a sequence of standard Pareto random variables with density (\ref{paretodf}) with $x_0=1$, therefore soft and hard tapering is defined as follows: if $b_n=n^\g$ and $\g>1/\a$, we have soft tapering, while $\g<1/\a$ gives us hard tapering. Clearly, only values $0<\a<2$ are interesting, and in the sequel we shall use this assumption without mentioning. Comparing with \cite{Chakrabarty}, were sums of tapered heavy-tailed random variables were considered, we consider more complicated objects - sums of linear processes with tapered innovations, therefore the picture is more complicated. Now the behavior of such sums is influenced not only by exponent $\a$ and tapering parameter $b_n$ but also  by a filter of a linear process. We know that for linear processes, depending on the filter we can have different asymptotic behavior of partial sums and four types of memory - positive, zero, negative, and strongly negative. We shall  combine the first three  of them with soft and hard tapering.
 %Also it is interesting what happens in the intermediate tapering when the limit (\ref{limtap}) is finite positive number (in \cite{Chakrabarty} intermediate truncation is considered). All these questions are left for the  future research.

%\section{Results}

We consider  partial sum processes of the form (\ref{linpr3}) and take $b_n=n^\g, \ \g>0$. Our goal is to investigate the
 asymptotic behavior of ${\bar S}_n(t)$ under soft and hard tapering.
 We do not seek for the most general conditions imposed on the filter of the linear process, therefore we assume that
\begin{equation}\label{cond1}
a_n\sim n^{-\b},
\end{equation}
where $\b >1/2$ (this condition ensures (\ref{cond})).  Changing the value of $\b$ and, in the case $\b>1$ making additional assumption $\sum_{i=0}^\infty a_i = 0$, we shall be able to consider all three cases of memory - positive, zero, and negative - for linear processes (with fixed tapering parameter $b$) under consideration. Our goal is to see what effect tapered innovations add in each case.

Intuitively it is clear that if the tapering parameter $b_n$ is growing rather slowly (hard tapering) the limiting process will be Gaussian, i.e. fBm $B_H$, while in the case of rapid growth of $b_n$ (soft tapering) limit
process will be stable. This intuition will be confirmed in two theorems formulated bellow.

Let us define $Z_n(t)=A_n^{-1}{\bar S}_n(t)$, where $A_n$ is a normalizing sequence for ${\bar S}_n(1)$.
Let $\{V_n\}\fdd \{V_0\}$ stand for the convergence of processes $V_n$ to a process $V_0$ in the sense of the finite-dimensional distributions (f.d.d.). In what follows letter $C$ will stand for   constants, dependent on various parameters, but not on $n$, and generally different in different places.

\begin{teo}\label{thm1} Suppose that the sum   (\ref{linpr3}) is formed by a linear process with a filter   (\ref{cond1})  and with  innovations   with parameter $0<\a<2$ and tapering parameter $b_n=n^\g, \ \g>0$, and one  the  following conditions holds:

%(i) if $1/2<\b<1$ and
\begin{equation}\label{cond5}
(i) \ if \ 1/2<\b<1 \ {\rm and} \ \g <\min \left (\frac{1}{\a}, \frac{2\b-1}{2-\a}\right );
\end{equation}
\begin{equation}\label{cond8}
(ii) \ if \ \b>1, \ \sum_{i=0}^\infty a_i \ne 0, \ {\rm and} \ \g <\min \left (\frac{1}{\a}, \frac{1}{2-\a}\right );
\end{equation}
\begin{equation}\label{condition11}
(iii)\ if \ 1<\b<3/2, \ \sum_{i=0}^\infty a_i = 0, \ {\rm and} \ \g <\min \left (\frac{2\b-1}{2-\a}, \frac{3-2\b}{\a}\right );
\end{equation}
then $A_n=Cn^H$, and, as $n\to \infty,$
$$\left \{Z_n(t) \right \} \fdd \{B_H( t)\},$$
where $B_H$ is fBm with parameter
\begin{equation}\label{Hcond}
H= \left \{\begin{array}{ll}
                     \frac{3}{2}-\b +\frac{\g(2-\a)}{2}  \ & {\rm in \ the \ cases} \ (i) \ {\rm and}\ (iii), \\
                      \frac{1}{2} +\frac{\g(2-\a)}{2}    \ & {\rm in\ the\ case} \ (ii).
                                                                       \end{array} \right.
\end{equation}

%\end{equation}
\end{teo}

\begin{remark}\label{rem1} Taking $\g=0$ in all three cases in the theorem (this means that we consider linear processes with a fixed $b$) we shall get well-known results for sums of linear processes with finite variance, see, for example, Propositions 3.3.1 and 4.4.1 in \cite{Giraitis}.
\end{remark}

\begin{remark}\label{rem2} Comparing the values of Hurst parameter $H$, given in all three cases in the theorem with those which we get taking $\g=0$, see Remark \ref{rem1}, we see that they had increased by the same value $\g(2-\a)/2$. Also it is worth to note that the boundary $1/\a$ for the hard tapering is achieved in cases (i) and (ii) and  for bigger values of $\a$: in case (i) for $\b^{-1}<\a<2$ and in the case (ii) for $1<\a<2$. This fact shows that if the Paretian part of distribution of innovations is heavier ($\a$ is smaller), the tapering sequence must grow slower. In the case (iii), since $(3-2\b)/\a<1/\a$,
the boundary $1/\a$ is not achieved. Let us note that the Brownian motion as limit process for sums (\ref{linpr3}) can be obtained only in the case (iii) with $\b=1+\g(2-\a)/2$.
\end{remark}

\begin{remark}\label{rem3} There is one more case of the so-called strongly negative memory, which, roughly speaking means that a filter of a linear process under consideration is such that normalization constants for $S_n(t)$  in (\ref{linpr5}) are bounded ($H=0$). Despite the fact that this case has some practical meaning (it corresponds to the over-differenced time series, see \cite{McElroy}), strongly negative memory is almost unexplored (in \cite{Paul20} there are examples of filters giving linear processes with such memory),  and we did not consider this case.
\end{remark}

Now we consider the soft tapering, this means that the tapering sequence $b_n$ is growing more rapidly and the Paretian part of distribution of innovations is playing more important role.  The idea is to compare ${\bar S}_n(t)$ with an appropriate sum of a linear process $Y_k=\sum_{j=0}^\infty a_j\eta_{k-j}, $ with Pareto innovations $\{\eta_i\}$, having the same exponent $\a$, as TP with the density (\ref{taperpd1}) and centered, if $\a>1$. Therefore, instead of the condition (\ref{cond}) now we must assume condition
\begin{equation}\label{cond7}
\sum_{j=0}^\infty |a_j|^\a <\infty,
\end{equation}
which gives us more narrow interval for parameter  $\b$. Moreover, we shall exclude the case $\a=1$, since in this case the condition (\ref{cond7}) does not ensure the a.s. convergence of the series $\sum_{j=0}^\infty a_j\eta_{k-j}$ and in the case $\a=1$ we need to require a little bit stronger condition $\sum_{j=0}^\infty |a_j|^{1-\e} <\infty$ with some $\e>0$. Although in our simple setting  (standard Pareto distribution) there is no principal difficulties dealing with the case $\a=1$, but in formulation of the result and in the proofs this case must be considered separately, therefore we decided to exclude this case.

To formulate our result we need some new notations. Let $L_\a(t), \ t\ge 0$, be an $\a$-stable process with independent stationary increments and characteristic function
\begin{equation}\label{Lalpha}
E\exp \{iuL_\a(t) \}=\exp \{-t|u|^\a(1-iD {\rm sign}u) \},\ \a \ne 1,
\end{equation}
with $D=-\tan \pi\a/2$. For $1/\a<\b<1+1/\a, \ \b\ne 1,$ let us define
\begin{equation}\label{LFSM}
Y(\a, \b, t)=\int_{-\infty}^t \left ((t-s)^{1-\b}-(-s)_+^{1-\b} \right )dL_\a(s), t\ge 0,
\end{equation}
and for $\b=1$ we set $Y(\a, 1, t)=L_\a(t).$
Introduced processes $L_\a$ and $Y(\a, \b, \cdot)$ are nothing else as $\a$-stable L{\' e}vy motion  and linear fractional stable motion, and it is possible to give different expressions via $\a$-stable random measures, see, for example, \cite{Samorod}. Since we shall cite one result from \cite{Astrauskas}, we use the definitions from this paper.

\begin{teo}\label{thm2}
Suppose that  the sum   (\ref{linpr3}) is formed by a linear process with a filter (\ref{cond1}),  satisfying (\ref{cond7}), and tapered innovations with parameters $1/\b<\a<2, \ \a \ne 1,$ and $b_n=n^\g, \ \g>0,$ and one of the following conditions holds:

(i) if  $1/\a<\b<1$ and  $\g>1/\a$;

(ii)  If $\b>1,  \ \sum_{i=0}^\infty a_i \ne 0$, and $\g>1/\a$;

(iii) if $\max (1, 1/\a)<\b<1+1/\a, \ \sum_{i=0}^\infty a_i = 0$,  and
\begin{equation}\label{cond12}
\g >\frac{1}{\a} +\frac{\b-1}{\a\b-1};
\end{equation}
 then, as $n\to \infty,$
\begin{equation}\label{convergfdd}
\left \{A_n^{-1}\sum_{j=1}^{[nt]}X_j^{(n)} \right \} \fdd \left \{\begin{array}{ll}
                                                                       \{Y(\a, \b, t)\}, \ & {\rm in \ the \ cases} \ (i)\ {\rm and}\ (iii), \\
                                                                       \{L_\a (t) \}, \ & {\rm in\ the\ case} \ (ii),
                                                                       \end{array} \right.
\end{equation}
where $Y(\a, \b, t)$ is defined in (\ref{LFSM}) and $L_\a(t)$ - in (\ref{Lalpha}), $A_n=Cn^H$ and $H=1/\a +1-\b$ in the cases (i) and (iii), and $H=1/\a$ in the case (ii).

\end{teo}

From \cite{Chakrabarty} we know that in the case of limit theorems for sums of truncated heavy-tailed random variables there is a simple dichotomy between hard and soft truncation (in our context it would be dichotomy between cases $\g<1/\a$ and $\g>1/\a$). In \cite{Chakrabarty} even the intermediate case, where the limit in (\ref{trunc}) is finite and positive number, is considered and limit law is obtained. In the case of linear processes with tapered innovations the situation is more complicated. At present there are even intervals of positive length of values of parameter $\g$, for which we do not know what is the limit law  for sums (\ref{linpr3}). In order to compare statements  (i)-(iii) from Theorem \ref{thm1} with corresponding  statements from Theorem \ref{thm2} we must take into account  that conditions on parameters $\a$ and $\b$ may be different due to the fact that different conditions (\ref{cond}) and (\ref{cond7}) are assumed in Theorems \ref{thm1} and \ref{thm2}, respectively. For example, if we want to compare hard and soft tapering in the case of negative memory, i.e., the cases (iii) in both theorems, we must consider the intersection of the intervals for the parameter $\b$, namely, intersection of $(1, 3/2)$ and $(\max (1, 1/\a), 1+1/\a).$ This intersection is non-empty for $2/3<\a<2$ and is equal
$$
(\max (1, 1/\a), 3/2)= \left \{\begin{array}{ll}
              (1/\a, 3/2), & {\rm if} \ 2/3<\a<1, \\
               (1, 3/2), & {\rm if} \ 1\le \a<2.
               \end{array} \right.
$$
Thus, only for these values of $\a$ and $\b$ we must compare two bounds from (\ref{condition11}) and (\ref{cond12}) for the parameter $\g$: $C_1(\a, \b):=\min \left ((2\b-1)(2-\a)^{-1}, \ (3-2\b)\a^{-1} \right )$ and $C_2(\a, \b):=\a^{-1}+  (\b-1)(\a\b-1)^{-1}$. It is not difficult to verify that for all values of $\a, \b$ under consideration there is a gap between these two constants, this means that for values of $\g$ in the interval $C_1(\a, \b)<\g < C_2(\a, \b)$ we do not know the limit process for sums (\ref{linpr3}). On the other hand, comparing statement (i) and (ii) and those $\a$ values, for which upper bounds in (\ref{cond5}) and (\ref{cond8}) are $1/\a$, we see that in these cases we have almost complete answer in a sense that only for $\g=1/\a$ the answer is not known.

\begin{remark}\label{rem4} It is possible to say that the choice of innovations with tapered distributions partially was motivated  by recent papers by L. Klebanov and his collaborators (see \cite{Klebanov1}, \cite{Klebanov2} ) where serious doubts about the usage of purely heavy-tailed distributions in finance is raised. In these papers it is demonstrated that some  distributions, like symmetrized gamma distribution, having exponential tails can explain many effects which usually are tried to explain by means of heavy-tailed distributions. Introducing tapered distributions we mentioned seismology, where tapered distributions occur quite naturally, in \cite{Chakrabarty} there are more examples of areas where tapered distributions are used to model real processes.
\end{remark}

\begin{remark}\label{rem5} Let us note that we do not touch the problem of statistical estimation of the parameters under consideration. If we observe tapered random variables itself, then the problem is easier, and in \cite{Chakrabarty} there are some statistical procedures to estimate tapering parameters $\a$ and $b$ are given, while in \cite{Aban} truncated Pareto distribution is considered and estimators for the parameters, present in (\ref{trunkPareto}) are given. It is interesting to note that in \cite{Aban} there are three log-log plots of empirical data from different fields, which are very similar to the Figure 1 from \cite{Kagan}, but in \cite{Aban} they all are considered as plots from truncated Pareto, while in \cite{Kagan} TP is used. In our case, if we observe only the values of a linear process $X_k(\a, b), 1\le k\le n,$ estimation of parameters $\a, \ b,$  also estimation  of filter parameter $\b$ seems to be a challenging problem for statisticians.
Another problem, which is left for the future research, is convergence of $Z_n(t)$ in some topology of the space $D[0,1]$.
\end{remark}

\section{Proofs}

 Let us denote
$\mu_r(b)=\int_1^\infty x^r f_b(x) {\dx},$
where $f_b(x)$ is from (\ref{taperpd1}). Easy calculations show that,  as $b\to \infty$,
\begin{equation}\label{taperpd2}
\mu_r(b)= \left \{\begin{array}{ll}
              \frac{r}{r- \a} b^{r-\a}(1+o(1)), & {\rm if} \ r>\a, \\
               \frac{\a}{\a-r} (1+o(1)), & {\rm if} \  r<\a, \\
              \a \ln b (1+o(1)), & {\rm if} \  r=\a.
               \end{array} \right.
\end{equation}
As usual, dealing with sums of values of a linear process, it is convenient  to write ${\bar S}_n(t)$ as infinite series of i.i.d random variables with weights:
\begin{equation}\label{sumrepres}
{\bar S}_n(t)=\sum_{k=1}^{[nt]}\sum_{j=0}^\infty a_j\xi^{(n)}_{k-j}=\sum_{j=-\infty}^{[nt]} d_{n,j,t}\xi^{(n)}_{j},
\end{equation}
where $d_{n,j,t}=\sum_{k=1}^{[nt]}a_{k-j}$,, for $j\le 0$ and $d_{n,j,t}=\sum_{k=j}^{[nt]}a_{k-j}$, for $j> 0,$ and $\xi^{(n)}_{k}$ are random variables $\xi_{k}$ with $b=b_n$.

Since from (\ref{sumrepres}) we have
\begin{equation}\label{sumvar}
{\rm Var}{\bar S}_n(1)=\sum_{j=-\infty}^n |d_{n,j}|^2 E(\xi^{(n)}_{1})^2,
\end{equation}
where $d_{n,j}:=d_{n,j, 1}$, and from (\ref{taperpd2}) we have the asymptotic of $E(\xi^{(n)}_{1})^2$, therefore we must investigate the asymptotic of $\sum_{j=-\infty}^n |d_{n,j}|^2$.

\begin{prop}\label{prop1}
Let the filter $\{a_n\}$ satisfy condition (\ref{cond1}), then
\begin{equation}\label{dnj1}
\sum_{j=-\infty}^n |d_{n,j}|^2 \sim \left \{\begin{array}{ll}
                                                                       C n^{3-2\b}, \ & {\rm in \ the \ cases} \ (i) {\rm and}\ (iii), \\
                                                                       C n, \ & {\rm in\ the\ case} \ (ii),
                                                                       \end{array} \right.
\end{equation}
where constants C in both relations depend only on $\b$ and are given  explicitly  in the proof.
\end{prop}

\begin{remark}\label{rem6} The quantity $\sum_{j=-\infty}^n |d_{n,j}|^2 $ is ${\rm Var}\sum_{j=1}^{n}X_j$, where $X_j, \ j\in \bz$ is a linear process with innovations , presenting standard white noise.
Asymptotic of this quantity is presented in \cite{Giraitis}, see Proposition 3.3.1 therein. It is interesting to note  that for the asymptotic in the case (iii) in \cite{Giraitis} a little bit stronger than condition (\ref{cond1}) is required, namely, $a_n =c_an^{-\b}(1+O(n^{-1})), \ c_a\ne 0,  \ 1<\b<3/2$. Probably stronger condition is needed, since in \cite{Giraitis} first the asymptotic of the  covariance of a linear process is derived from asymptotic of $a_n$,  then the asymptotic of variance of $\sum_{j=1}^{n}X_j$ is obtained, see Propositions 3.2.1 and 3.3.1. In our approach we use the asymptotic of $a_n$ directly to establish the asymptotic of the variance of  the sum $\sum_{j=1}^{n}X_j$.
\end{remark}

{\it Proof of Proposition \ref{prop1}}. We divide the sum under the considerations into two parts
\begin{equation}\label{dnj3}
\sum_{j=-\infty}^n |d_{n,j}|^2=V_1 +V_2:=\sum_{j=-\infty}^0 |d_{n,j}|^2 + \sum_{j=1}^n |d_{n,j}|^2.
\end{equation}
We can write $a_n=n^{-\b}(1+\d(n))$, where $\d(n)\to 0$, as $n\to \infty$, then
$$
V_1=\sum_{j=0}^\infty \left (\sum_{k=1}^n a_{k+j} \right )^2=\sum_{j=0}^\infty \sum_{k=1}^n\sum_{l=1}^n(k+j)^{-\b}(l+j)^{-\b}\left (1+\d_1(k, l, j) \right ),
$$
where $\d_1(k, l, j)\to 0$, as $j\to \infty$, uniformly over all $k\ge 0, l\ge 0$. Changing sums into integrals and making change of variables we arrive at the following quantity:
$$
\int_0^\infty\int_1^n\int_1^n (x+w)^{-\b}(y+w)^{-\b}\left (1+\d_1(x, y, w) \right )dx dy dw=
$$
$$
n^{3-2\b}\int_0^\infty\int_{1/n}^1\int_{1/n}^1 (u+z)^{-\b}(v+z)^{-\b}\left (1+\d_1(nu, nv, nz) \right )du dv dz.
$$
Now we prove that in the cases (i) and (iii)
\begin{equation}\label{dnj4}
V_1\sim n^{3-2\b}\int_0^\infty\left (\int_{0}^1(u+z)^{-\b} du\right)^2  dz.
\end{equation}
First we show that for $1/2<\b<3/2, \b\ne 1$ the multiple integral in (\ref{dnj4}) is finite. Writing the integral in (\ref{dnj4}) as a sum
\begin{equation}\label{dnj5}
\int_0^1\left (\int_{0}^1(u+z)^{-\b} du\right)^2 dz  + \int_1^\infty\left (\int_{0}^1(u+z)^{-\b} du\right)^2 dz,
\end{equation}
consider the case $1/2<\b<1$. In the first integral we use the bound $\int_{0}^1(u+z)^{-\b} du<C$ and prove the finiteness of the first integral, while the second can be written as
$$
\int_1^\infty \left ( (1+z)^{1-\b}-(z)^{1-\b}\right )^2 dz=\int_1^\infty z^{-2\b}\left (\left (\frac{z} {1+z}\right )^{\b}(1+z)-z\right )^2 dz.
$$
Is not difficult to show that the function $\left (\left (\frac{z} {1+z}\right )^{\b}(1+z)-z\right )$ is bounded, therefore the second integral in (\ref{dnj5}) is also bounded,  thus, the multiple integral in (\ref{dnj4}) in  the case $1/2<\b<1$ is finite. In the case  $1<\b<3/2$  for the first integral we simply use the fact that $\int_0^1z^{2(1-\b)}dz<C $, while in order to show the finiteness of the second integral we can use the same expression as above.

To finish the proof of (\ref{dnj4}) we must  show that the integral
$$
\int_0^\infty\int_{1/n}^1\int_{1/n}^1 (u+z)^{-\b}(v+z)^{-\b}\d_1(nu, nv, nz)du dv dz
$$
can be made arbitrary small for sufficiently large $n$. To this aim we divide the integral $\int_0^\infty$ with respect to $z$ into two integrals $\int_0^\e\dots dz$ and $\int_\e^\infty\dots dz$, where $\e>0$ is small but fixed number. In the second integral $\int_\e^\infty\dots dz$ we have $\d_1(nu, nv, nz)\to 0$, as $n\to \infty$, uniformly with respect to all possible values of $u, v$, therefore  the integral
$$
\int_\e^\infty \int_{1/n}^1\int_{1/n}^1 (u+z)^{-\b}(v+z)^{-\b}\d_1(nu, nv, nz)du dv dz
$$
can be made arbitrary small for sufficiently large $n$. Finally, let us note that
$$
\int_0^\e \int_{1/n}^1\int_{1/n}^1 (u+z)^{-\b}(v+z)^{-\b}\d_1(nu, nv, nz)du dv dz \to 0, \quad {\rm as } \quad \e\to 0,
$$
since in the case $1/2<\b<1$ this integral is of the order $C\int_0^\e dz$ and in the case $1<\b<3/2$ it is of the order $C\int_0^\e z^{2(1-\b)}dz$ and $-1<2(1-\b)<0$. Thus we have proved (\ref{dnj4}).
In the case $1<\b$ and $\sum_{i=0}^\infty a_i \ne 0$ we shall show that
\begin{equation}\label{dnj9}
\frac{V_1}{n} \to 0.
\end{equation}
For this it is sufficient to estimate $|a_i|\le Ci^{-\b}$ and using integral criterion to get the bounds:
\begin{eqnarray*}
\frac{V_1}{n} & \le &\frac{C}{n}\left (\sum_{j=0}^n \left (\sum_{k=1}^n a_{k+j} \right )^2+\sum_{j=n+1}^\infty \left (\sum_{k=1}^n a_{k+j} \right )^2\right )\\
     & \le &\frac{C}{n}\sum_{k=1}^n k^{1-\b} +\frac{C}{n}\sum_{j=n+1}^\infty nj^{-2\b}\le Cn^{1-\b}.
\end{eqnarray*}

Investigation of the term $V_2$ is more simple. Let us consider the case $1/2<\b<1$. Repeating the same steps as in investigation of $V_1$ we have
$$
V_2=\int_0^n\int_0^{n-w}\int_0^{n-w} x^{-\b}y^{-\b}\left (1+\d_1(x, y, w) \right )dx dy dw=
$$
$$
n^{3-2\b}\int_0^1\int_0^{1-w}\int_0^{1-w} u^{-\b}v^{-\b}\left (1+\d_1(nu, nv, nz) \right )du dv dz.
$$
Since $\b<1$ we easily arrive at the relation
\begin{equation}\label{dnj6}
V_2\sim n^{3-2\b}\int_0^1\left (\int_0^{1-z} u^{-\b}du\right )^2 dz.
\end{equation}
Consider the case $1<\b<3/2$ and $\sum_{i=0}^\infty a_i=0$. Using this assumption and the same procedure as in earlier considerations we arrive at the following relation
$$
V_2=\sum_{i=1}^n\left ( \sum_{k=n-j+1}^\infty a_k\right )^2=n^{3-2\b}\int_0^1\left (\left (\int_{1-z}^\infty u^{-\b}du \right)^2+\d_2(n, 1-z) \right )dz,
$$
where $\d_2(n, 1-z)\to 0$, for $n\to \infty$, uniformly over $1-z>\e$ with a fixed small $\e$. The main term in the asymptotic relation is $\int_0^1\left (\int_{1-z}^\infty u^{-\b}du \right)^2dz$ which is finite due to assumption $-1<2(1-\b)<0$. Since the outer integral is over $0\le z\le 1$, to deal with the remainder term, as earlier we must divide integral into two integrals over $0\le z <1-\e$ and $1-\e\le z \le 1$, and finally we get
\begin{equation}\label{dnj7}
V_2\sim n^{3-2\b}\int_0^1\left (\int_{1-z}^\infty u^{-\b}du \right)^2dz.
\end{equation}
The case $1<\b$ and $\sum_{i=0}^\infty a_i\ne 0$ is the most simple, since using the Toeplitz lemma (see, for example, \cite{Loeve}, p. 250) we have
$$
\frac{V_2}{n}=\frac{1}{n}\sum_{k=0}^{n-1}\left ( \sum_{j=0}^k a_j\right )^2 \to \left ( \sum_{j=0}^\infty a_j\right )^2,
$$
thus, we have
\begin{equation}\label{dnj8}
V_2\sim n\left (\sum_{j=0}^\infty a_j \right)^2.
\end{equation}

It remains to collect the relations (\ref{dnj3}), (\ref{dnj4}), (\ref{dnj9}), (\ref{dnj6})-(\ref{dnj8}) and we get (\ref{dnj1}).
\halmos

{\it Proof of Theorem \ref{thm1}}.
We start with the proof of  Theorem \ref{thm1} in the case (i). From (\ref{sumvar}), taking into account (\ref{taperpd2}), it is easy to get, as \ $n\to \infty$,
\begin{equation}\label{varbn}
{\rm Var}{\bar S}_n(1)=\sum_{j=-\infty}^n |d_{n,j}|^2 E(\xi^{(n)}_{1})^2\sim C b_n^{2-\a}n^{3-2\b}=Cn^{2H},
\end{equation}
where $H$ is from (\ref{Hcond}).
Let us note that the condition
\begin{equation}\label{cond3}
0<\g<\frac{2\b-1}{2-\a},
\end{equation}
 ensures that  $0<H<1$. It is not difficult to show that under this condition we have
\begin{equation}\label{covzn}
{\rm Var}Z_n(t) \to t^{2H}, \  \ {\rm Cov}(Z_n(t), Z_n(s)) \to \frac{1}{2}\left (|s|^{2H}+|t|^{2H}-|s-t|^{2H} \right ),
\end{equation}
and this is the covariation function of a fBm $B_H$. Now we prove the convergence of f.d.d. of $Z_n$ to corresponding f.d.d. of $B_H$.
As usual, the convergence of f.d.d. is proved with the help of Cram{\' e}r-Wold device. Since this step became standard in similar questions, therefore to simplify the writing, instead of considering linear combination $d_1Z_n(t_1)+ \dots +d_k Z_n(t_k)$ we shall take $k=1, d_1=1$ and even $t_1=1$.
To prove that the distribution of $Z_n(1)$ converges to the standard normal law it is sufficient, according to Lyapunov CLT, to show that, as $n\to \infty,$
 \begin{equation}\label{Liapunov}
L(3,n)=\frac{\sum_{j=-\infty}^n |d_{n,j}|^3}{\left (\sum_{j=-\infty}^n d_{n,j}^2\right )^{3/2}}\frac{E|\xi^{(n)}_{1}|^3}{\left (E(\xi^{(n)}_{1})^2\right )^{3/2}} \to 0.
\end{equation}
Estimating $\sum_{j=-\infty}^n |d_{n,j}|^3\le \max_j |d_{n,j}|\sum_{j=-\infty}^n |d_{n,j}|^2$ and remembering the expression of $d_{n, j}$, we easily get
\begin{equation}\label{Liapunov1}
\frac{\sum_{j=-\infty}^n |d_{n,j}|^3}{\left (\sum_{j=-\infty}^n d_{n,j}^2\right )^{3/2}}\le \frac{\max_j |d_{n,j}|}{\left (\sum_{j=-\infty}^n |d_{n,j}|^2\right )^{1/2}}\le C\frac{n^{1-\b}}{n^{3/2-\b}}=Cn^{-1/2}.
\end{equation}
Taking into account (\ref{taperpd2}) and our choice $b_n=n^\g$, we have
\begin{equation}\label{Liapunov2}
\frac{E|\xi^{(n)}_{1}|^3}{\left (E(\xi^{(n)}_{1})^2\right )^{3/2}}\le C\frac{b_n^{3-\a}}{b_n^{3(2-\a)/2}}=Cn^{\g\a/2}.
\end{equation}
If the condition
\begin{equation}\label{cond4}
\g <\frac{1}{\a}
\end{equation}
holds, then (\ref{Liapunov1}) and (\ref{Liapunov2}) prove (\ref{Liapunov}). Combining (\ref{cond3}) and (\ref{cond4}) we have   bound (\ref{cond5}) for $\gamma$, and the case (i) is proved.

Considering the case (ii), the proof goes along the same lines as the proof of the case (i), therefore we shall provide  the changes only. Instead of (\ref{varbn}) now we have
%\begin{equation}\label{varBbn}
$${\rm Var}{\bar S}_n(1)=\sum_{j=-\infty}^n |d_{n,j}|^2 E(\xi^{(n)}_{1})^2\sim C b_n^{2-\a}n=Cn^{2H},$$
%\end{equation}
where  $H$ is from (\ref{Hcond}).
It is easy to see that condition
$0<\g<(2-\a)^{-1}$ ensures that $0<H<1$. Instead of (\ref{Liapunov1})  now we have
%\begin{equation}\label{Liapunov1B}
$$\frac{\sum_{j=-\infty}^n |d_{n,j}|^3}{\left (\sum_{j=-\infty}^n d_{n,j}^2\right )^{3/2}}\le \frac{\max_j |d_{n,j}|}{\left (\sum_{j=-\infty}^n |d_{n,j}|^2\right )^{1/2}}\le Cn^{-1/2}$$
%\end{equation}
and (\ref{Liapunov2}) remains unchanged, therefore under condition (\ref{cond4}) we get (\ref{Liapunov}) The rest steps in the proof are the same as in the proof of the case (i).

It remains to prove the case (iii). From Proposition \ref{prop1} we have that in the case (iii) the asymptotic of $\sum_{j=-\infty}^n |d_{n,j}|^2 $ is the same as in the case (i) (only we must stress that if in the case (i) parameter $1/2<\b<1$, now we have  $1<\b<3/2$), therefore we get the same relation (\ref{varbn}) with $H$ from (\ref{Hcond}).

  Requiring $0<H<1$ we get the same condition  (\ref{cond3}) and in the same way we get (\ref{covzn}). It remains to prove the convergence of f.d.d., that is, to show (\ref{Liapunov}). Writing the equality
$$
\max_{-\infty<j\le n}|d_{n,j}|=\max (\max_{-\infty<j\le 0}|d_{n,j}|, \ \max_{1\le j\le n}|d_{n,j}|),
$$
 using the expression of $d_{n, j}$ and relation $|\sum_{k=0}^{n-j}a_{k}| =|\sum_{k=n-j+1}^{\infty}a_{k}|$, we easily get
$$
\max_{-\infty<j\le n}|d_{n,j}|\le \sum_{k=1}^{\infty}|a_{k}|.
$$
Now we have
\begin{equation}\label{Liapunov1C}
\frac{\sum_{j=-\infty}^n |d_{n,j}|^3}{\left (\sum_{j=-\infty}^n d_{n,j}^2\right )^{3/2}}\le \frac{\max_j |d_{n,j}|}{\left (\sum_{j=-\infty}^n |d_{n,j}|^2\right )^{1/2}}\le Cn^{-3/2+\b}.
\end{equation}
Since we have (\ref{Liapunov2}), combining it and (\ref{Liapunov1C}) we get (\ref{Liapunov}), if condition
\begin{equation}\label{cond4C}
\g <\frac{1}{\a}(3-2\b)
\end{equation}
holds. Combining (\ref{cond3}) and (\ref{cond4C}) we have (\ref{condition11}). The case (iii) is proved.
\halmos

{\it Proof of Theorem \ref{thm2}}.
Let us introduce the standard Pareto random variable $\theta$ with the density
$$
g(x)=\a x^{-\a-1},  \ x\ge 1,
$$
and let $R$ be standard exponential random variable with the  density function $e^{-x}$, for $x\ge 0$, independent of $\theta$. Then it is easy to see that a random variable $\zeta$ with the density (\ref{taperpd1}) can be written as
$$
\zeta=\theta\ind{[\theta<b]}+(b+R)\ind{[\theta\ge b]}.
$$
Let us denote $\xi=\zeta-E\zeta$
and $\eta=\theta -E\theta$, if $1<\a <2$ and $\eta=\theta$, if $0<\a < 1$.
For $\kappa<\a$ it is not difficult to obtain the following estimate
$$
E|\eta -\xi|^\kappa\le Cb^{-(\a-\kappa)},
$$
where $C$ depends on $\a, \kappa$.

 Let $\{\eta_i,\ i\in \bz\}$ and $\{\xi_i, \ i\in \bz\}$  be  sequences of i.i.d. copies of $\eta$ and $\xi$, respectively,  defined on the same probability space as  $\theta$, $\zeta$, and  in such a way that for all pairs  $\eta_i, \xi_i$ the above written inequality holds, namely,
\begin{equation}\label{compkappa}
E|\eta_i -\xi_i|^\kappa\le Cb^{-(\a-\kappa)}.
\end{equation}
 Let us consider a linear processes
$Y_k=\sum_{j=0}^\infty a_j\eta_{k-j}, \ k\in \bz,$
and corresponding partial sum process $V_n(t)=\sum_{k=1}^{[nt]}Y_k, \ 0\le t\le 1$. We start with the case (i).
 Due to conditions  imposed on the filter $\{a_i\}$, for the process $V_n$ we can apply statement (ii) from Theorem 1 in \cite{Astrauskas} to get  that
\begin{equation}\label{convergfddZ}
\left \{V_n(t)(Cn^{H})^{-1} \right \} \fdd \{Y(\a, \b, t)\},
\end{equation}
where $Y(\a, \b, t) $ is defined in (\ref{LFSM}) and $H=1/\a +1-\b$. Since the exponent of normalization constants is of the form $H=1/\a +\d$ with $\d=1-\b>0$,
we  have positive memory case (again using terminology from \cite{Paul20}).    The  skewness parameter $D,$ reflecting balance between tails of distribution of innovations, has such value due to the fact that Pareto distribution is one-sided. In order to prove (\ref{convergfdd}) we shall show that, for  any $\e>0$ and any fixed $t$,
\begin{equation}\label{convsergP}
\PP\{n^{-H}|V_n(t)-\sum_{k=1}^{[nt]}X_k^{(n)}|>\e \} \pc 0.
\end{equation}
The relation (\ref{convsergP})  will follow from the relation
\begin{equation}\label{convsergE}
n^{-H\kappa}E|V_n(t)-\sum_{k=1}^{[nt]}X_k^{(n)}|^\kappa \to 0, \ \ {\rm for \ some }\ \ \kappa>0.
\end{equation}
Since we intend to use (\ref{compkappa}), we  take $1<1/\b<\kappa<\a$, the lower bound $1/\b<\kappa$ is used to ensure that $\sum_{i=0}^\infty |a_i|^\kappa<\infty$.
Denoting by $\xi^{(n)}_{k}$ the random variable $\xi_{k}$ with tapering parameter $b_n,$ we have
$$
V_n(t)-\sum_{k=1}^{[nt]}X_k^{(n)}=\sum_{k=1}^{[nt]}\sum_{j=0}^\infty a_j(\eta_{k-j}-\xi^{(n)}_{k-j}).
$$
For the  sum over $k$ we use the rough estimate
$$
E|V_n(t)-\sum_{k=1}^{[nt]}X_k^{(n)}|^\kappa\le [nt]^{\kappa-1}\sum_{k=1}^{[nt]}E|\sum_{j=0}^\infty a_j(\eta_{k-j}-\xi^{(n)}_{k-j})|^\kappa,
$$
while for the infinite sum over $j$ we use 2.6.20 result from \cite{Petrov} (we recall that $E(\eta_i -\xi_i)=0$ and $1<\kappa<2$) and then
(\ref{compkappa}). We get
$$
n^{-H\kappa}E|V_n(t)-\sum_{j=1}^{[nt]}X_j^{(n)}|^\kappa \le Cn^{-H\kappa+\kappa-\g(\a-\kappa)}.
$$
 Substituting the value of $H$ and requiring that the exponent at $n$ would be negative, we get
$$
\g>g(\kappa):=\frac{\kappa(\a \b-1)}{\a(\a-\kappa)}.
$$
Since the value of $\kappa$ can be chosen from interval $1/\b<\kappa<\a$, we need to find $\inf_{1/\b<\kappa<\a}g(\kappa)$. It is easy to verify that
$g(1/\b)=1/\a$  and  $\lim_{\kappa\to \a}g(\kappa)=\infty$,
and the derivative of the function $g$ is positive in the interval under consideration, therefore, we get
$$
\inf_{1/\b<\kappa<\a}g(\kappa)=\frac{1}{\a}.
$$
We cannot take $\kappa=1/\b$ (due to the fact that $\sum_{i=0}^\infty |a_i|^{1/\b} =\infty$), but we can take the value of $\kappa$ arbitrary close to $1/\b$. Thus,  if $\g>1/\a$, then we can choose value of $\kappa$ arbitrary close to $1/\b$, for which  (\ref{convsergE}) holds, and the case (i) is proved.

The proof of the case (ii) goes along the same lines as the proof of the case (i), again  we shall provide  only the changes.
We consider the introduced sums $V_n(t)=\sum_{j=1}^{[nt]}Y_j, \ 0\le t\le 1$, and now, applying statement (i) from Theorem 1 in \cite{Astrauskas}, we get  that
\begin{equation}\label{convergfddZB}
\left \{n^{-1/\a}V_n(t) \right \} \fdd \{L_\a(t)\}.
\end{equation}
Having (\ref{convergfddZB}) we must prove (\ref{convsergP}), only with $H=1/\a$. Let us consider the case $1/\b<\a \le 1$, then we choose $1/\b<\kappa<\a \le 1$ and we get
$$
E|V_n(t)-\sum_{j=1}^{[nt]}X_j^{(n)}|^\kappa\le \sum_{k=1}^{[nt]}E|X_k-Y_k|^\kappa)\le nt \sum_{i}a_i^\kappa b_n^{-(\a-\kappa)}\le Cn^{1-\g(\a-\kappa)}.
$$
For this choice of $\kappa$ and if $\g>1/\a$, we get
\begin{equation}\label{convsergEB}
n^{-\kappa/\a}E|V_n(t)-\sum_{j=1}^{[nt]}X_j^{(n)}|^\kappa\le Cn^{1-\g(\a-\kappa)-\kappa/\a} \to 0.
\end{equation}
If $1<\a<2$ and $\b>1$ we can take $\kappa=1$ and in the same way, if $\g>1/\a$, we get (\ref{convsergEB}). Therefore, we have (\ref{convsergP}) (with $H=1/\a$), this relation, together with (\ref{convergfddZB}), proves the case (ii).

It remains to prove (iii). Now we have negative memory for the process $V_n$, therefore, applying statement (iii) from Theorem 1 in \cite{Astrauskas}, we get (\ref{convergfddZ}) with the same $H$, only now $H=1/\a+\d$ with $\d=1-\b<0$. Now, in order to prove (\ref{convsergP}), we must get (\ref{convsergE}). To this aim we choose $1/\b<\kappa<\min (\a, 1)$ and, applying the same inequalities as in the proof of the case (ii), we get
$$
n^{-H\kappa}E|V_n(t)-\sum_{j=1}^{[nt]}X_j^{(n)}|^\kappa\le Cn^{1-\g(\a-\kappa)-H\kappa}.
$$
Substituting the value of $H$ into the expression of the exponent in the last estimation and requiring that this exponent would be negative we shall get
$$
\g>h(\kappa):=\frac{\a-\left (1+\a(1-\b) \right )\kappa}{\a (\a-\kappa)}.
$$
It is easy to verify that
$$
\inf_{1/\b<\kappa<\min(\a, 1)}h(\kappa)=h(1/\b)=\frac{1}{\a}+\frac{\b-1}{\a\b-1}.
$$
 Therefore,  having condition (\ref{cond12}), we get (\ref{convsergE}) and, together with  (\ref{convergfddZ}), we prove (\ref{convergfdd}). The theorem is proved.
\halmos

\end{document}